\documentclass[10pt]{article}
\usepackage{latexsym}

\begin{document}

\begin{center}
{\bf\large  Note on the number of proper colorings of a graph}

\bigskip
{\large Martin Klazar}\\
{\em
Department of Applied Mathematics\\
and\\ 
Institute for Theoretical Computer Science (ITI)\\
Faculty of Mathematics and Physics, Charles University\\
Malostransk\'e n\'am\v est\'\i\ 25, 118 00 Praha\\
Czech Republic\\
{\tt klazar at kam.mff.cuni.cz}
}
\end{center}

\bigskip\noindent
Using a counting argument called Tur\'an sieve (with motivation  in number theory), 
Liu and Murty proved 
in \cite[Theorem 4]{liu_murt} an inequality on the number of proper
vertex colorings, by $\lambda$ colors, of a simple graph $G=(V,E)$ with $v=|V|$ vertices and $e=|E|$
edges:
$$
\#(\mbox{proper $\lambda$-colorings of $G$})\le\lambda^v\cdot\frac{\lambda-1}{e}.
$$
This note presents a more direct combinatorial argument 
that gives a slightly stronger inequality.

\bigskip
\noindent
{\bf Theorem.} If $G=(\{1,2,\dots,v\},E)$ is a simple graph with $v\ge 1$ vertices and 
$e=|E|\ge 0$ edges and $\lambda\ge 1$ is the number of colors, then
$$
|C^p|:=\#(\mbox{proper $\lambda$-colorings of $G$})\le
\lambda^v\cdot\frac{\lambda-1}{e+\lambda-1}.
$$

\smallskip\noindent
{\em Proof.} For $e=0$ the inequality holds and so does for $\lambda=1$ (for $\lambda=1$ and 
$e=0$ we interpret $0/0$ as $1$), and we assume that $\lambda\ge 2$ 
and $e\ge 1$. Let $L$, $|L|=\lambda$, be the set of colors, $C$ be the set of 
colorings of $V$ by colors from $L$, and $C^p\subset C$ be the colorings which are 
proper (i.e., have no monochromatic edge). We shall construct an injection
$$
I:\;\{(h,f):\;h\in E, f\in C^p\}\to\{(c,g):\;c\in L, g\in C\backslash C^p\}
$$
with the property that for every coloring $g\in C\backslash C^p$ there are at most
$\lambda-1$ colors $c$ such that $(c,g)$ is in the image ${\rm Im}(I)$ of $I$. Then 
we must have
$$
e|C^p|\le (\lambda-1)(|C|-|C^p|)=(\lambda-1)(\lambda^v-|C^p|)
$$
and the inequality follows.

\medskip

For every subset $X\subset V$, fix a spanning forest $F_X$ of the graph $G_X$ 
induced by $G$ on $X$ (i.e., fix a 
spanning tree in every component of $G_X$). Now for $(h,f)\in E\times C^p$, where 
$h=\{u,w\}$ with $u<w$ and $f(u)=d\ne f(w)=c$, consider the set $Y\subset V$
of vertices having in $f$ color $d$ or $c$. In the spanning forest
$F_Y$, $u$ and $w$ lie in the same component $K$ and $K$ contains a unique path
$P$ joining $u$ and $w$. Recolor the vertices of $K$ by $c$ and $d$ so that 
all vertices of $P$ are colored with $d$ and the only monochromatic edges in $E(K)$
(and hence in $E(F_Y)$) are those 
of $P$; there is exactly one such coloring of $K$ (since $K$ is a tree). 
Define $g\in C$ as given by this recoloring on $K$ and as 
coinciding with $f$ elsewhere. Set $I((h,f))=(c,g)$; recall that $c$ is the color 
lost on $P$ by the recoloring. Clearly, $(c,g)\in L\times C\backslash C^p$. Note that the 
only {\em bad} colors in $g$---colors appearing
on some edge in $E$ monochromatic in $g$---are $d$ (because of the edges in $P$) and possibly $c$ 
(the recoloring may have created some edges monochromatic in the color $c$ but these 
must lie in $E\backslash E(F_Y)$). The membership $(c,g)\in {\rm Im}(I)$ imposes on 
$(c,g)$ further restrictions which we discuss in a moment.

\medskip

We have to verify the injectivity of $I$ and the property of its image.
Let $(c,g)\in (L\times C\backslash C^p)\cap {\rm Im}(I)$ be given. We describe how to 
reconstruct uniquely $(h,f)$ from $(c,g)=I((h,f))$. First, set $B'=\{c\}\cup B$ where 
$B$ is the set of colors bad in $g$. As we mentioned, necessarily $|B'|=2$. Thus $Y$ 
is uniquely reconstructed as $Y=g^{-1}(B')$. The color $d\in B'$, $d\ne c$, is the only 
color that is bad in $g$ with respect to the edges in $E(F_Y)$ and, moreover, the 
monochromatic edges in $E(F_Y)$ form a path $P$ with endvertices $u<w$ that must be adjacent
(this is forced by $(c,g)\in {\rm Im}(I)$). Thus the 
edge $h$ is uniquely reconstructed as $\{u,w\}$. If we modify 
$g$ on the component $K$ of $F_Y$ containing $P$ so that $K$ is properly colored
with $c$ and $d$ but $u$ retains its color $d$ (such a coloring of $K$ is unique), 
the resulting coloring is necessarily proper with respect to all edges in $E$ 
(else $(c,g)$ could not be in ${\rm Im}(I)$) and $f$ must be equal to it. Hence we have reconstructed 
$(h,f)$. This shows that $I$ is an injection.  

\medskip

Finally, let $(c,g)\in {\rm Im}(I)$ and $B$ be the set of bad colors in $g$. 
If $|B|=1$, then $c\not\in B$ and for fixed $g$ we have at most $\lambda-1$ possibilities for $c$. 
If $|B|=2$, then $c\in B$ and is uniquely
determined (must be distinct from the color in $B$ that is bad in $g$ with respect to the edges 
in $E(F_Y)$ where $Y=g^{-1}(B)$) and there are again at most $\lambda-1$ possibilities for $c$ 
when $g$ is fixed. 
\hfill{$\Box$}

\bigskip\noindent
{\bf Closing remarks.} Stronger result was obtained long ago by Lazebnik in \cite{laze}: 
$$
\#(\mbox{proper $\lambda$-colorings of $G$})\le\lambda^v\cdot A
$$
where 
$$
A=\min\left(
\bigg(1-\frac{1}{\lambda}\bigg)^{\lceil\sqrt{2e+1/4}-1/2\rceil},\ 
1-\frac{e}{\lambda}+\frac{{e\choose 2}}{\lambda^2},\ 
\frac{\lambda-1}{e+\lambda-1}
\right).
$$
The last third term in the minimum is identical with the one obtained here. I could not check 
\cite{laze} to see if the above argument based on an injective mapping, 
which still may be of some interest, is subsumed in the proofs of \cite{laze}. 
For further strengthening and references on the 
problem to bound the number of proper colorings, see Byer \cite{byer}.

\end{document}